\documentclass[a4paper,12pt]{article}%Autor: Becherer, Bernhardt, Gapeev

\usepackage{amsmath}%many things (e.g. equation in theorem command, align, gather ...)
\usepackage{amssymb}%R, C, Q with double-bar, {\setminus} otherwise strange spacing depending on sets like $R_+\setminus D$ and $[0,T[\,\setminus D$
\usepackage{amsthm}%\newtheorem command
\usepackage{bbm}%indicator function (1 with double-bar, missing in amssymb)
\usepackage[
    backend=biber, % using helper program for biblatex (often not needed)
    articlein=false, % removes 'In:' from articles, only available for ext-styles
    dashed=false, % removes '-' and displays author names when two articles are written by the same authors
    giveninits=true,% abbreviated first names
    maxcitenames=2, % et al when there are more than two names in citing
    natbib=true, % macros from incompatible package natbib like 'citep' or 'citet'
    style=ext-authoryear-ibid, % yields author (year) in Reference list
    uniquename=init,% technical commands to avoid warning (ext-authoryear-ibid uses full name while giveninits reduces full name)
    ]{biblatex}%bibliography creator and reference list changes
\usepackage[margin=2.5cm]{geometry}%changing page layout, 3cm edge, symmetric 
\usepackage[low-sup]{subdepth}%turning off the automated height adjustments for sub and super scripts (now everything on the same height)

\usepackage{todonotes}%comments in paper
\usepackage{xcolor}%colors

\allowdisplaybreaks%breaks up math environments over pages otherwise only on the same page
\citetrackerfalse %no ``ibid'' same source
\newtheoremstyle{NoItalic}%how theorems and definition appear
    {1em}%space above
    {1em}%space below
    {}%body font
    {}%indent
    {\bfseries}%theorem head font
    {.}%punctuation after theorem head
    {.5em}%space after theorem head
    {\thmname{#1}\thmnumber{ #2}\thmnote{ (#3)}}%theorem head
\theoremstyle{NoItalic}

\newtheorem{defi}{Definition}[section]
\newtheorem{ex}[defi]{Example}
\newtheorem{lem}[defi]{Lemma}
\newtheorem{rmk}[defi]{Remark}
\newtheorem{thm}[defi]{Theorem}

\makeatletter
\renewenvironment{proof}[1][\proofname]{\par
\pushQED{\qed}%
\normalfont \topsep6\p@\@plus6\p@\relax
\trivlist
\item\relax
{\itshape
#1\,\normalfont\@addpunct{:}}\hspace\labelsep\ignorespaces
}{%
\popQED\endtrivlist\@endpefalse
}
\makeatother

\DefineBibliographyStrings{english}{and={\&}} % replacing 'and' with '&'
\DeclareFieldFormat[article,incollection]{title}{#1} % removes quotation marks for title in article
\DeclareFieldFormat[article,periodical]{number}{\mkbibparens{#1}} % volume numbers are in round brackets
\DeclareNameAlias{sortname}{family-given} % surname-firstname in bibliography
\begin{document}

\title{On Watanabe's characterisation and change of intensity \`{a} la Girsanov for Cox processes}
\author{Dirk Becherer\footnote{Humboldt-Universität zu Berlin, {\texttt{dirk.becherer@hu-berlin.de}}}, \hspace{0.15em} Thomas Bernhardt\footnote{University of Manchester, {\texttt{thomas.bernhardt@manchester.ac.uk}}} \hspace{0.15em} and Pavel Gapeev\footnote{London School of Economics and Political Science, {\texttt{p.v.gapeev@lse.ac.uk}}}}
%\date{}
\maketitle

\begin{abstract}
\noindent
We discuss the equivalence of definitions for conditional Poisson processes, Cox processes, and stochastic intensities of point processes on the real line. A basic proof shows that Watanabe's characterisation of conditional Poisson processes in terms of local martingales is necessary and sufficient. 
% In particular, we show that conditional Poisson processes are Cox processes. 
%D: steht nun schon in erster zeile  %T: agree
Additionally, we consider conditions enabling the measure change method \`a la Girsanov to alter the intensity of Cox processes to a desired new target intensity, e.g.\ for the probability reference approach in filtering. Such a measure change exists if a corresponding stochastic exponential is a proper martingale. We show that this holds
% under the minimal condition, that 
% D:minimal conditions weglassen? %T: agree
if the new locally integrable target intensity is the product of the original intensity and another non-negative process. 
\\[1em]
\\[1em]
\textbf{Acknowledgements:} D.Becherer acknowledges funding by the Deutsche Forschungsgemeinschaft (DFG, German Research Foundation) - CRC/TRR 388 "Rough Analysis, Stochastic Dynamics and Related Fields" - Project ID 516748464.
\end{abstract}

%---------------------------
\section{Introduction}\label{section:Intro}
%---------------------------

We consider a simple point process $N$ on the positive half-line in a complete and filtered
%D: check questions in todo on last page %T: done 
probability space $(\Omega,\mathcal{F},(\mathcal{F}_t)_{t\geq0},\mathbb{P})$, which means that there are stopping times $(T_n)_{n=0}^\infty$ such that
\begin{equation}\label{eq:N}
    \begin{gathered}
        N_t=\sum_{n\in\mathbb{N}}\mathbbm{1}_{\{T_n\leq t\}}\quad\text{for $\,t\geq0$},
        \\\mathbb{P} \big[ T_0=0 \big] = 1 \quad\text{and}\quad\mathbb{P} \big[ T_{n+1}<\infty\;\Rightarrow\;T_n<T_{n+1} \big] = 1 \quad\text{for $\,n\in\mathbb{N}_0$},
    \end{gathered}
\end{equation}
and a non-negative $\,\mathcal{B}([0,\infty))\otimes\mathcal{F}$-measurable process $X$ that is locally integrable, i.e.\
\begin{equation}\label{eq:X}
    \mathbb{P}\Big[\int_0^t\!X_s\;\mathrm{d}s<\infty\quad\forall\;t\geq0\Big]=1.
\end{equation}
First, we show that the following three definitions are equivalent, provided that $X$ is an $\mathcal{F}_0$-measurable process.
%D: (in the sense that the defining properties are equivalent), T: did not understand what the note in brackets does and removed it  
This can be viewed as a generalisation of the famous characterisation by \textcite[remark to Theorem 2.3]{Watanabe1964} for the Poisson process, which corresponds to the special case where the intensity process $\,X\equiv1\,$ is constant.

\begin{defi}\label{defi:conditional Poisson}
    $N$ is a conditional Poisson process with $\mathcal{F}_0$-measurable intensity $X$ if for all $\,u\in\mathbb{R}\,$ and $\,t\geq r\geq0$
    \[\mathbb{E}\Big[\exp\Big(iu(N_t-N_r)\Big)\,\Big|\,\mathcal{F}_r\Big]=\exp\Big((\mathrm{e}^{iu}-1)\int_r^tX_s\,\mathrm{d}s\Big).\]
\end{defi}

\begin{defi}\label{defi:Cox}
    $N$ is a Cox process with $\mathcal{F}_0$-measurable intensity $X$ if for all finitely many real-valued constants $(u_j)_{j=1}^k$ and bounded disjoint $\,\mathcal{B}([0,\infty))$-measurable sets $(A_j)_{j=1}^k$
    \[\mathbb{E}\Big[\exp\Big(i\sum_{j=1}^ku_j\int_{A_j}\!dN_s\Big)\,\Big|\,\mathcal{F}_0\Big]=\exp\Big(\sum_{j=1}^k(\mathrm{e}^{iu_j}-1)\int_{A_j}\!X_s\,\mathrm{d}s\Big).\]
\end{defi}

\begin{defi}\label{defi:mloc-intensity}
    $N$ has a stochastic $(\mathcal{F}_t)$-intensity $X$ if the  compensated point process $\,\widetilde{N}_t=N_t-\int_0^tX_s\,\mathrm{d}s$, $\,t\geq0$, is a local $(\mathcal{F}_t)$-martingale.
\end{defi}

\noindent Second, we show that the stochastic exponential
\begin{gather}\label{eq:Z}
    Z_t=\exp\Big( -\int_0^t (Y_s-1) X_s \, \mathrm{d}s\Big)\, 
    \prod_{s\in]0,t]}\!\!\big( 1+(Y_s-1) \, \Delta N_s \big),\quad t\geq0,
\end{gather}
is a martingale with respect to the underlying filtration $(\mathcal{F}_t)_{t\geq0}$ when $N$ is a conditional Poisson process with $\mathcal{F}_0$-measurable intensity $X$ and $Y$ is a non-negative $\,\mathcal{B}([0,\infty))\otimes\mathcal{F}_0$-measurable process such that
\[
\mathbb{P}\Big[\int_0^tY_sX_s\,\mathrm{d}s<\infty\quad\forall\;t\geq0 \Big]=1.
\]

%D: changed jump process to counting process; check for consistence in remainder of paper. jump pr is too general. pure jump levy process may be non monotone, non counting, have exploding counts of jumps....} %T: agree, by the way I could not find any other times when 'jump process(es)' was mentioned 
Counting processes, such as the Poisson process, are commonly used to model the counts of recurring events. The literature (e.g.\ \textcite[Chapter~13]{LastPenrose2017}, or \textcite[Chapters~1 and 5]{Bremaud2020}) describes such processes by defining an intensity that determines the average frequency of events per unit of measurement.
%D: per unit of time measurement; T: not necessarily time (can be anything); similar stochastic evolution instead of stochastic time evolution
The possibly stochastic evolution of the intensity corresponds to fluctuations in the (locally expected) frequency of point events. 
%Often, the intensity itself is a stochastic process. 
Typically, one uses one of the three Definitions \ref{defi:conditional Poisson}, \ref{defi:Cox} or \ref{defi:mloc-intensity} to describe counting processes, and these are (indeed) equivalent without further conditions for counting processes in the setting described above with an $\mathcal{F}_0$ measurable intensity $X$. We note for comparison that the occurrence of jumps in such processes does not affect the future dynamics of the intensity. In contrast, that would be different for point processes with dimming or self-exciting features, like Hawkes processes.  

%In this note, we look at conditional Poisson processes on the real line. 
%D: steht bereit weiter oben %T: agree 

To validate both the first and second statements above, we employ only basic arguments and rather direct computations, which are based on the observation that conditional expectations of event stopping times $(T_n)_{n=0}^\infty$ can be calculated in a fairly basic closed form, in the sense of Lemma~\ref{lem:E[Y_T_n+1]=}. 
%
% D hier ergaenzt, es ist nun im sinne das Fazit wieder etwas mehr wie in meiner erste Aenderung, hoffentlich nun klarer lesbar. Siehe auch im naechsten abchnitt Ergaenzung not "citable" mit "easily accessible"
%
We wish to emphasize that the present paper serves to provide elementary proofs for the aforementioned two fundamental statements on Cox processes. We believe that such basic results deserve elementary proofs. Notwithstanding that, we are of course not suggesting, that those statements could not alternatively be derived also from more advanced theory, e.g., by employing less elementary results involving so-called characteristics from general semimartingale theory (see \textcite[Chapter~II]{JS} and also note Remark~\ref{rmk:deterministic-variation}).  

%%%%%%%%%%%
%T: changed this sentence a couple of times ... I hope it reads better now ... what do you think, Dirk?
%%%%%%%%%%%%

%T: JS assumes Definition 1.3 together with X F_0 measurable, in particular, that reference misses the point; I have put CR together with KallsenKarbe

%We wish to emphasize that the present short note serves to provide an elementary proof for the aforementioned fundamental statements on Cox processes, which does not require to employ results from general semimartingale theory (see, e.g., \textcite[Section~II.6]{JS}, the more recent article\footnote{We have been made aware of this reference after the first arXiv version of this note was published.} by \textcite[Theorem~4.1]{CR}, or \textcite{KallsenKarbe2010} as described in Remark~\ref{rmk:deterministic-variation}).

\hfill

%-------------------------
\noindent\textbf{Equivalence of definitions (see Theorem~\ref{thm:Watanabe}):} The equivalence of the three Definitions \ref{defi:conditional Poisson}, \ref{defi:Cox} and \ref{defi:mloc-intensity} under the assumption that $X$ is an $\mathcal{F}_0$-measurable process might well be folklore among experts. %Yet, despite conditional Poisson processes and Cox processes being a well-studied topic, it has been surprisingly hard for us to find any citable, easily accessible reference for the equivalence of the three definitions above. Indeed, when both process classes have appeared in publications, they are often treated separately, or their equivalence is shown under additional assumptions.
Yet, despite conditional Poisson processes and Cox processes being a well-studied topic, it has been surprising to us that we have not been able to find a readily citable publication, which directly states 
the equivalence of the three definitions above 
and moreover presents the result in an explicitly accessible way. Let us note  that the identity of the processes $X$ in the three distinct characterizing properties forms a part of the equivalence to be established.
% . Indeed, when both process classes have appeared in publications, they are often treated separately, or their equivalence is shown under additional assumptions.   

Definitions~\ref{defi:Cox} and \ref{defi:mloc-intensity} are equivalent, provided that $X$ is an $\mathcal{F}_0$-measurable process (see \cite[Example 5.1.5 and Exercise 5.9.3]{Bremaud2020}). Clearly, Def.~\ref{defi:Cox} implies Def.~\ref{defi:conditional Poisson}. However, it is non-trivial that Def.~\ref{defi:conditional Poisson} implies Def.~\ref{defi:Cox} or \ref{defi:mloc-intensity}.

A first guess that a monotone class (or Dynkin) argument may show the implication from Def.~\ref{defi:conditional Poisson} to Def.~\ref{defi:Cox} turns out to fail, at least without additional arguments.
% D: folgendes wirklich ausfuehren? falls ja (gut), bitte praeziser erklaeren 
% T: Tried to do my best to explain the issue briefly; if it is does not work, lets rather drop it (but at least mention that we thought about it) 
%It is true that finite unions of the right-closed intervals $\,(r,t]\,$ for $\,t\geq r\geq0\,$ form an intersection-stable generator of $\,\mathcal{B}([0,\infty))$ and a successive application of the tower property to an increasingly ordered list of intervals using Def.~\ref{defi:conditional Poisson} yields Def.~\ref{defi:Cox}. However, the crux is how to generalise to countably many disjoint sets. Here, the lack of an increasing order between the sets is an obstacle to combining the conditional expectations corresponding to different disjoint sets into a single conditional expectation. 

\citet[Theorem 4.1]{Grigelionis1975} considered the implication Def.~\ref{defi:conditional Poisson} to Def.~\ref{defi:mloc-intensity}.
%D: check elsewhere in paper: do not ref. to plain label numbers only, without def eqn etc % T: done
He states that both definitions are equivalent, provided that $X$ is an $\mathcal{F}_0$-measurable process. He states ``necessity of the theorem's conditions follows from [his] definition 4.1'' in his proof on page 447. However, his Definition 4.1 does not include a martingale property at first sight. It is worth mentioning that \citeauthor{Grigelionis1975} appears to make the standing assumption that $\,t\mapsto\int_0^tX_s\,\mathrm{d}s\,$ always refers to the compensator on page 445 and, as such, assumes the martingale property. But crucially, he did not show a connection between the process $X$ as in Def.~\ref{defi:conditional Poisson} and $X$ as in Def.~\ref{defi:mloc-intensity}.

Due to the one-to-one correspondence between probability distributions and characteristic functions, the equation in Def.~\ref{defi:conditional Poisson} is equivalent to 
\[
    \mathbb{P}\big[N_t-N_r=n\,\big|\,\mathcal{F}_r\big]=\frac{(\int_r^tX_s\,\mathrm{d}s)^n}{n!}\mathrm{e}^{\textstyle-\int_r^tX_s\,\mathrm{d}s}\quad\text{for all $\,n\in\mathbb{N}_0$}.
\]
As both sides of the above equation are positive, 
monotone convergence implies that Def.~\ref{defi:conditional Poisson} (see also \cite[Exercise 5.9.9]{Bremaud2020}) yields
\begin{equation}\label{eq:not-enough-?}
    \mathbb{E}\big[ N_t-N_r\,\big|\,\mathcal{F}_r\big]=\int_r^t X_s\,\mathrm{d}s.
\end{equation}
Equation~(\ref{eq:not-enough-?}) clearly implies Def.~\ref{defi:mloc-intensity} if either side of the equation is integrable. For example, \textcite[Chapter~II.2]{Bremaud1981} assumes $\,\mathbb{E}[N_t]<\infty\,$ and \textcite[Proposition~18]{Grandell1991} assumes $\,\mathbb{E}[\int_0^tX_s\,\mathrm{d}s]<\infty$. 
However, without assuming the integrability of either side in Equation~(\ref{eq:not-enough-?}), one needs to introduce a localising sequence $(T_n)_{n=0}^\infty$ to make either side integrable. Indeed, as mentioned in \textcite[Remark 5.1.4]{Bremaud2020}, the correct condition that implies Def.~\ref{defi:mloc-intensity} is 
\begin{equation}\label{eq:now-enough}
    \mathbb{E}\big[N_{t\wedge T_n}-N_{r\wedge T_n}\,\big| \,\mathcal{F}_r\big]=\mathbb{E} \Big[\int_{r\wedge T_n}^{t\wedge T_n}\!X_s\,\mathrm{d}s\,\Big| \,\mathcal{F}_r\Big].
\end{equation} 
But Equation~(\ref{eq:not-enough-?}) does not imply Equation~(\ref{eq:now-enough}) in an obvious way because optional sampling does not apply here, as that would require a martingale in the first place. An alternative shortcut to this end may be to use a sequence of $\mathcal{F}_0$-measurable stopping times instead of $(T_n)_{n=0}^\infty$; we will look at this alternative route in Remark~\ref{rmk:F0-mb-time}.

\hfill

\noindent\textbf{Change of measures (see Theorem~\ref{thm:Change-of-Intensity}):} The question of whether $Z$ as given in Equation~(\ref{eq:Z}) is not just a local martingale but a martingale determines when we can use $Z$ as a Radon-Nikodym density to define a measure change under which the intensity of the conditional Poisson process is changed from $X$ to $YX$, see \textcite[Chapter~VI]{Bremaud1981}. For example, this is a crucial assumption in the reference probability approach to filtering for point processes. Here, the goal is to extract information about $X$ from observation $N$. That might be easy or difficult depending on the complexity of the joint distribution of $X$ and $N$. Now, the reference probability approach simplifies the joint distribution by introducing a measure change that makes $X$ independent of $N$. Then, additional properties, like convenient marginal distributions, allow us to compute conditional expectations of $X$ given $N$ using the Kallianpur--Striebel formula (see \textcite{Zakai1969} or \textcite[Chapter~3]{BainCrisan2009} for the approach in a diffusion setting). We exemplify the reference probability approach in Example~\ref{ex:MSc}. 
%D: you wanted to cut out the filter examples from arxiv ver1 in earlier version - are they introduced again now? if not, consider to refer to arxiv ver1 for the detailed examples} %T: Yes, there were two examples, I included the more interesting one (but not both; showcasing with one example is good enough) 

Whether a stochastic exponential is a martingale or only a local martingale is of mathematical interest on its own. We refer the reader to \textcite{MijatovicUrusov2012} or \textcite{Ruf2013} and the references therein for the literature on this topic in a diffusion setting. The known conditions in the point process setting typically impose bounds on the integrands involved (and work irrespective of any class of point processes). Examples can be found in \textcite[Section~3.3]{BoVaWo1975}, \textcite[Theorem 2.4]{HansenSokol2015}, or \textcite[Remark 5.5.2 and subsequent examples]{Bremaud2020}. We will explain in   Remark~\ref{rmk:deterministic-variation} also a notable result by \textcite{KallsenKarbe2010}
%\textcite[Theorem 3.1, Proposition 3.12]{KallsenKarbe2010}, 
which implies our claim under slightly stronger assumptions in the context of general semimartingale theory, and also a recent generalization by \textcite{CR} that removed the additional assumptions. 

Let us emphasise once more that the present note serves to provide elementary proofs for the aforementioned  basic statements on Cox processes, without employing less elementary techniques or more advanced results like from the general theory of semimartingales (cf.\ \cite{JS}).

%---------------------------
\section{Statements and Proofs}
%---------------------------

The following lemma provides the main basic tool for our analysis. It shows how the conditional expectations at event times of conditional Poisson processes can be written down in closed form. This will allow us to calculate any relevant conditional expectation to check martingale properties in later theorems.

\begin{lem}[Conditional Density]\label{lem:E[Y_T_n+1]=}
    Let $N$ be a conditional Poisson process with $\mathcal{F}_0$-measurable intensity $X$, let $Y$ be a $\,\mathcal{B}([0,\infty])\otimes\mathcal{F}_r$-measurable non-negative process for given $\,r\geq0$, and $\,n\in\mathbb{N}_0$. Then
    \[\mathbb{E} \big[ Y_{T_{n+1}}\, \big| \,\mathcal{F}_r \big] =
        \begin{cases}
            \int_r^\infty Y_t\Psi_t^{n,r}\mathrm{d}t+Y_\infty{\displaystyle\Psi_\infty^{n,r}}&\mbox{on $\,\{N_r\leq n\}$}
        \\
            \hfil Y_{T_{n+1}}&\mbox{on $\,\{N_r>n\}$}
        \end{cases}
    \]
    with conditional density $\,\Psi^{n,r}\colon[0,\infty]\times\Omega\to[0,\infty)\,$ given by
    \[\Psi_t^{n,r}=
        \begin{cases}
            \displaystyle X_t\,\mathrm{e}^{\textstyle-\int_r^tX_s\,\mathrm{d}s} \, \frac{1}{(n-N_r)!} \, \Big( \int_r^tX_s\,\mathrm{d}s \Big)^{n-N_r}\;\mathbbm{1}_{\{N_r\leq n\}}&\mbox{if $\,t<\infty$}
        \\
            \displaystyle\hfil\mathrm{e}^{\textstyle-\int_r^\infty X_s\,\mathrm{d}s}
    			\sum_{k=0}^{n-N_r}\frac{1}{k!} \, \Big(\int_r^\infty X_s\,\mathrm{d}s \Big)^k\;\mathbbm{1}_{\{N_r\leq n\}}&\mbox{if $\,t=\infty$}
       \end{cases}
    \]
    using the conventions $\,0^0=1\,$ and $\,\mathrm{e}^{-\infty}\infty=0$.
\begin{proof}
    We prove the statement for the special case $\,Y_s=\mathbbm{1}_{\{s>t\}}\,$ for given finite $\,t\geq r\geq0$. The general result then follows from the measurability and a monotone class argument.

    Because $N$ is a conditional Poisson process with intensity $X$, we have
    \begin{align*}
        \mathbb{P} \big[ T_{n+1}>t\, \big| \,\mathcal{F}_r \big]&=\mathbb{P}\big[N_t\leq n\,\big|\,\mathcal{F}_r\big]=\mathbb{P}\big[N_t-N_r\leq n-N_r\,\big|\,\mathcal{F}_r\big]
        \\&=\sum_{k=0}^{n-N_r}\mathbb{P}\big[N_t-N_r=k\,\big|\,\mathcal{F}_r\big]
        \\&=\mathrm{e}^{\textstyle-\int_r^tX_s\,\mathrm{d}s}
			\sum_{k=0}^{n-N_r}\frac{1}{k!} \, \Big( \int_r^tX_s\,\mathrm{d}s \Big)^k 
            \; \mathbbm{1}_{\{N_r\leq n\}}.
   \end{align*}
    The above process is absolutely continuous with respect to $t$ because the sums, products, and compositions (with increasing absolutely continuous functions) of absolutely continuous functions are again absolutely continuous. In particular, there is a derivative $\Psi^{n,r}$ such that for all finite $\,v\geq t\geq r$
    \begin{gather}
        \label{eq:int_Psi=}\int_t^v\Psi_s^{n,r}\,\mathrm{d}s=      \mathbb{P}\big[T_{n+1}>t\,\big|\,\mathcal{F}_r\big]-
        \mathbb{P}\big[T_{n+1}>v\,\big|\,\mathcal{F}_r\big]\quad\mbox{and}
        \\\label{eq:Psi=}\Psi_t^{n,r}=X_t\,\mathrm{e}^{\textstyle-\int_r^tX_s\,\mathrm{d}s}
        \, \frac{1}{(n-N_r)!} \Big(\int_r^tX_s\,\mathrm{d}s \Big)^{n-N_r}\; 
        \mathbbm{1}_{\{N_r\leq n\}}.
    \end{gather}
    Furthermore, consider $\,\Psi_\infty^{n,r}=\mathbb{P}[T_{n+1}=\infty\,|\,\mathcal{F}_r]\,$, then the monotone convergence theorem and $\,x^k\mathrm{e}^{-x}\rightarrow0\,$ when $\,x\uparrow\infty\,$ for any $\,k\in\mathbb{N}_0$ shows that
    \begin{equation}\label{eq:Psi_infty=}
        \Psi_\infty^{n,r}=\lim_{v\rightarrow\infty}
        \mathbb{P}\big[T_{n+1}>v\,\big|\,\mathcal{F}_r\big]=\mathrm{e}^{\textstyle-\int_r^\infty X_s\,\mathrm{d}s}\sum_{k=0}^{n-N_r}\frac{1}{k!} \Big(\int_r^\infty X_s\,\mathrm{d}s \Big)^k \; \mathbbm{1}_{\{N_r\leq n\}}.
    \end{equation}
    Combining equations \eqref{eq:int_Psi=}, \eqref{eq:Psi=}, \eqref{eq:Psi_infty=} yields the required result
    \[\mathbb{P}\big[ T_{n+1}>t\, \big| \,\mathcal{F}_r\big]=
    \int_t^\infty \Psi_s^{n,r}\,\mathrm{d}s+\Psi_\infty^{n,r}=\int_r^\infty Y_s\,\Psi_s^{n,r}\,\mathrm{d}s+Y_\infty\,\Psi_\infty^{n,r}.\]
\end{proof}
\end{lem}

Next, we show that Def.~\ref{defi:conditional Poisson} implies Def.~\ref{defi:mloc-intensity}.
%from the introduction. 
Let us formulate the result as a generalisation of \textcite{Watanabe1964}'s seminal characterization of the Poisson process. 

\begin{thm}[Characterisation \`{a} la \textcite{Watanabe1964} and \textcite{Grigelionis1975}]\label{thm:Watanabe}
    Let $N$ be a simple point process corresponding to $(T_n)_{n=0}^\infty$ and let $X$ be an $\mathcal{F}_0$-measurable process. Then, the following statements are equivalent
    \begin{itemize}
        \item[(i)] $N$ is a conditional Poisson process with intensity $X$ (see Def.~\ref{defi:conditional Poisson}),
        \item[(ii)]$\widetilde{N}_t=N_t-\int_0^tX_s\,\mathrm{d}s\,$ for $\,t\geq0\,$ is a local martingale.
    \end{itemize}
\begin{proof}
    The implication (ii) $\Rightarrow$ (i) is well-known, see for example \textcite[Chapter II.2]{Bremaud1981}. Thus, we only show the implication (i) $\Rightarrow$ (ii).

    Using the notation of Lemma~\ref{lem:E[Y_T_n+1]=}, we have for finite $\,t\geq r\geq0\,$ that
    \begin{align}
        \mathbb{E} \Big[\int_{r\wedge T_n}^{t\wedge T_n}\!X_s\,\mathrm{d}s\, \Big| \,\mathcal{F}_r \Big] \nonumber
        &=\mathbb{E} \Big[ \mathbbm{1}_{\{T_n>t\}}\int_r^t X_s\,\mathrm{d}s+\mathbbm{1}_{\{t\geq T_n>r\}}\int_r^{T_n}\!X_s\,\mathrm{d}s\, \Big| \,\mathcal{F}_r \Big]
        \\&=\mathbb{P}\big[T_n>t\,\big|\,\mathcal{F}_r \big] \, \int_r^tX_s\,\mathrm{d}s+\int_r^t\Psi_u^{n-1,r}\!\int_r^uX_s\,\mathrm{d}s\;\mathrm{d}u.
        \label{eq:E_int_Tn_X=}
    \end{align}
    Using the special form of $\Psi$ in particular $\,\Psi_u^{n-1,r}\int_r^uX_s\,\mathrm{d}s=(n-N_r)\Psi_u^{n,r}$, it holds that
    \begin{equation}\label{eq:int_Psi_int_X=}
        \int_r^t\Psi_u^{n-1,r}\int_r^uX_s\,\mathrm{d}s\;\mathrm{d}u = (n-N_r)\,\int_r^t\Psi_u^{n,r}\,\mathrm{d}u=(n-N_r)\;\mathbb{P}\big[ t\geq T_{n+1}>r\,\big|\,\mathcal{F}_r\big].
    \end{equation}
    By the definition of a conditional Poisson process, we have
    \begin{align}
        &\mathbb{P}\big[T_n>t\,\big|\,\mathcal{F}_r\big] \, \int_r^tX_s\,\mathrm{d}s
        =\mathbb{P}\big[N_t<n\,\big|\,\mathcal{F}_r\,\big]\int_r^tX_s\,\mathrm{d}s\nonumber
        \\&=\mathbb{P}\big[N_t-N_r<n-N_r\,\big|\,\mathcal{F}_r\big]\,\int_r^tX_s\,\mathrm{d}s
        =\mathrm{e}^{\textstyle-\int_r^tX_s\,\mathrm{d}s}\sum_{k=0}^{n-N_r} k\; \frac{1}{k!} \Big( \int_r^tX_s\,\mathrm{d}s\Big)^k \nonumber
        \\&=\sum_{k=0}^{n-N_r}k\;\mathbb{P}\big[N_t-N_r=k\,\big|\,\mathcal{F}_r\big]
        =\mathbb{E}\big[(N_t-N_r) \, \mathbbm{1}_{\{T_{n+1}>t\}}\,\big|\,\mathcal{F}_r\big].\label{eq:P_int_X=}
    \end{align}
    Moreover, the following equation holds,
    \begin{equation}
    \label{eq:N_Tn=}
    N_{t\wedge T_n}=N_{r\wedge T_n}+(n-N_r)\,\mathbbm{1}_{\{t\geq T_{n+1}>r\}}+(N_t-N_r)\,\mathbbm{1}_{\{T_{n+1}>t\}}.
    \end{equation}
    Now, combining equations \eqref{eq:E_int_Tn_X=}, \eqref{eq:int_Psi_int_X=}, \eqref{eq:P_int_X=}, \eqref{eq:N_Tn=} and noting that all following terms are quasi-integrable (possibly $-\infty$ but never $\infty$) or bounded yields
    \begin{align*}
        &\mathbb{E}\Big[N_{t\wedge T_n}-\int_0^{t\wedge T_n}X_s\,\mathrm{d}s \, \Big|
        \,\mathcal{F}_r\Big]=N_{r\wedge T_n}-\int_0^{r\wedge T_n}X_s\,\mathrm{d}s
        \\&+\mathbb{E}\Big[(n-N_r)\,\mathbbm{1}_{\{t\geq T_{n+1}>r\}}
        +(N_t-N_r)\,\mathbbm{1}_{\{T_{n+1}>t\}}-\int_{r\wedge T_n}^{t\wedge T_n}X_s\,\mathrm{d}s\,\Big|\,\mathcal{F}_r\Big]\\
        &=N_{r\wedge T_n}-\int_0^{r\wedge T_n}X_s\,\mathrm{d}s.
    \end{align*}
    Thus, the process $\,t\mapsto N_t-\int_0^tX_s\,\mathrm{d}s\,$ stopped with $T_n$ fulfils the martingale property; it is integrable (which follows from the special case $\,r=0$); and it is obviously adapted to $(\mathcal{F}_t)_{t\geq0}$; i.e.\ it is a martingale.

    Note that $\,T_n\uparrow\infty\,$ for $\,n\uparrow\infty\,$ because Def.~\ref{defi:conditional Poisson} implicitly implies that increments of $N$ are finite, i.e.\ $\mathbb{P}[N_t<\infty]\,$ for all $\,t\geq0$. In particular, $(T_n)_{n=0}^\infty$ is a suitable localising sequence.
\end{proof}
\end{thm}

\begin{rmk}\label{rmk:each-jump}
    Applying the optional stopping theorem to Lemma~\ref{lem:E[Y_T_n+1]=} shows that for any $\mathcal{B}([0,\infty])\otimes\mathcal{F}_{T_n}$-measurable non-negative process $W$ for given $\,n\in\mathbb{N}_0$
    \begin{equation}\label{eq:EW_Tn+1=}
        \mathbb{E} \big[W_{T_{n+1}}\,\big|\,\mathcal{F}_{T_n}\big]=
        \int_{T_n}^\infty W_s\, \Phi_s^n \,\mathrm{d}s+W_\infty \, \Phi_\infty^n
    \end{equation}
    with $\,\Phi^n\colon[0,\infty]\times\Omega\to[0,\infty)\,$ given by
    \[\Phi_s^n=
        \begin{cases}
            X_s\,\mathrm{e}^{\textstyle-\int_{T_n}^s X_r\,\mathrm{d}r}&\mbox{if $\,s<\infty$}
        \\
            \hfil\mathrm{e}^{\textstyle-\int_{T_n}^\infty X_r\,\mathrm{d}r}&\mbox{if $\,s=\infty$}
       \end{cases}.
    \]
    
    Now, we can also show the implication (i) $\Rightarrow$ (ii) from Theorem~\ref{thm:Watanabe} by applying \eqref{eq:EW_Tn+1=} to general point process theory. \textcite{ChouMeyer1975} (see also \textcite[Theorem 5.2.2]{Bremaud2020})
    developed the idea of computing the compensator of the entire point process from the compensators of every single jump. More precisely, in our context,
    \[
      \lambda_t=\sum_{n=0}^\infty\frac{\Phi_t^n}{1-\int_{T_n}^t\Phi^n_s\,\mathrm{d}s}
      \,\mathbbm{1}_{\{T_{n+1}>t\geq T_n\}}\quad\mbox{for $\,t\geq0$}
    \]
    has the property that $\,t\mapsto N_t-\int_0^t\lambda_s\,\mathrm{d}s\,$ is a local martingale. Plugging $\Phi^n$ in the above equation reveals $\,\lambda_t=X_t$.
\end{rmk}

\begin{rmk}\label{rmk:F0-mb-time}
    Another way to show the implication (i) $\Rightarrow$ (ii) from Theorem~\ref{thm:Watanabe} is to use a sequence of $\mathcal{F}_0$-measurable stopping times instead of $(T_n)_{n=0}^\infty$. Let
    \begin{align*}
        S_n&=\inf\Big\{t\geq0\,\Big|\,\int_0^tX_s\,\mathrm{d}s\geq n\Big\},
        \\\hat{S}_n&=\lfloor n\times S_n\rfloor/n\quad\text{for $\,n\in\mathbb{N}$},
    \end{align*}
    in which $\lfloor\cdot\rfloor$ denotes the floor function. Because $\hat{S}_n$ takes only a countable number of values, Equation~(\ref{eq:not-enough-?}) implies
    \begin{align*}
        \mathbb{E}\big[N_{t\wedge\hat{S}_n}-N_{r\wedge\hat{S}_n}\,\big|\,\mathcal{F}_r\big]&=\mathbb{E}\Big[\sum_{k=0}^\infty (N_{t\wedge\hat{S}_n}-N_{r\wedge\hat{S}_n})\mathbbm{1}_{\hat{S}_n=k/n}\,\Big|\,\mathcal{F}_r\Big]
        \\&=\sum_{k=0}^\infty\mathbb{E}\big[N_{t\wedge(k/n)}-N_{r\wedge(k/n)}\,\big|\,\mathcal{F}_r\big]\mathbbm{1}_{\hat{S}_n=k/n}
        \\&=\sum_{k=0}^\infty\int_{r\wedge(k/n)}^{t\wedge(k/n)}\!X_s\,\mathrm{d}s\;\mathbbm{1}_{\hat{S}_n=k/n}
        \\&=\int_{r\wedge\hat{S}_n}^{t\wedge\hat{S}_n}\!X_s\,\mathrm{d}s\quad\text{for $\,t\geq r\geq 0$}.
    \end{align*}
    By definition, $\,\hat{S}_n\leq S_n$, in particular, $\,\int_0^{t\wedge\hat{S}_n}\!X_s\,\mathrm{d}s\,$ is bounded, which means that we can subtract or add it to any well-defined expectation. Hence,
    \begin{gather*}
        \mathbb{E}\big[N_{t\wedge\hat{S}_n}-\int_0^{t\wedge\hat{S}_n}\!X_s\,\mathrm{d}s-\big(N_{r\wedge\hat{S}_n}-\int_0^{r\wedge\hat{S}_n}\!X_s\,\mathrm{d}s\big)\,\big|\,\mathcal{F}_r\big]=0 \quad\text{and}
        \\\mathbb{E}\big[N_{t\wedge\hat{S}_n}-\int_0^{t\wedge\hat{S}_n}\!X_s\,\mathrm{d}s\,\big|\,\mathcal{F}_0\big]=0\quad\text{for $\,r=0\,$ (implying integrability)}.
    \end{gather*}
    Thus, $\,t\mapsto\widetilde{N}_{t\wedge\hat{S}_n}=N_{t\wedge\hat{S}_n}-\int_0^{t\wedge\hat{S}_n}\!X_s\,\mathrm{d}s\,$ is a martingale. 

    $\hat{S}_n$ would not necessarily be a stopping time if $X$ was only adapted to $(\mathcal{F}_t)_{t\geq0}$, because $\hat{S}_n$ stops $\,\int_0^tX_s\,\mathrm{d}s\,$ before $\,\int_0^tX_s\,\mathrm{d}s\,$ reveals that it reaches $n$. However, $\mathcal{F}_0$ belongs to all the sigma algebras $\mathcal{F}_t$. Thus, any positive $\mathcal{F}_0$-measurable random variable is a stopping time, including $\hat{S}_n$.

    Furthermore, $\,S_n-1/n\leq\hat{S}_n\,$ and $\,S_n\uparrow\infty\,$ for $\,n\rightarrow\infty\,$ by definition. So does $\,\hat{S}_n\uparrow\infty\,$ for $\,n\rightarrow\infty$. Combining all the above statements yields that $(\hat{S}_n)_{n=1}^\infty$ is a localising sequence for the process $\,t\mapsto\widetilde{N}_{t}$. In particular, $\,t\mapsto\widetilde{N}_{t}\,$ is a local martingale with respect to filtration $(\mathcal{F}_t)_{t\geq0}$. 
\end{rmk}

%------------ 2nd Issue
Similarly to the proof of Theorem~\ref{thm:Watanabe}, we can check the martingality of other processes involving the conditional Poisson process $N$ using Lemma~\ref{lem:E[Y_T_n+1]=}. Next, we look at stochastic exponentials involving $\widetilde{N}$, which are essential for change of measure results in the context of Girsanov transformations.

\begin{thm}[Changing intensity by change of measure]\label{thm:Change-of-Intensity}
    Let $N$ be a conditional Poisson process with $\mathcal{F}_0$-measurable intensity $X$, let $Y$ be a non-negative $\,\mathcal{B}([0,\infty))\otimes\mathcal{F}_0$-measurable process such that
    \[
    \mathbb{P}\Big[\int_0^t Y_sX_s \,\mathrm{d}s<\infty\quad\forall\;t\geq0
    \Big]=1.\]
    Consider the stochastic exponential $\,Z=\mathcal{E}(\int_0^\cdot Y_s-1\,\mathrm{d}\widetilde{N}_s)\,$ with the compensated process $\widetilde{N}$ from Theorem~\ref{thm:Watanabe}(ii), that is,
    \[
    Z_t=\exp\Big( -\int_0^t (Y_s-1) X_s \, \mathrm{d}s\Big)\, 
    \prod_{s\in]0,t]}\!\!\big( 1+(Y_s-1) \, \Delta N_s \big) \quad \mbox{for $\,t\geq0$}\]
    with the convention $\,\prod_\varnothing=1$. 
    Then, the following statements are true:
    \begin{itemize}
        \item[(i)]$Z$ is a non-negative right-continuous martingale with $\,Z_0=1$.
        \item[(ii)]
            The probability measure induced by $Z$ and a bounded $\mathcal{F}_0$-measurable random variable $\,S\geq0$, more precisely,
            \[\mathbb{Q}(A)=\mathbb{E}\big[\mathbbm{1}_A\,Z_S\big]\quad\mbox{for $\,A\in\mathcal{F}$},\]
            is such that $\,\mathbb{Q}=\mathbb{P}\,$ on $\mathcal{F}_0$ and the process $\widetilde{N}_{t\wedge S}^{\mathbb{Q}}=N_{t\wedge S}-\int_0^{t\wedge S} Y_sX_s\,\mathrm{d}s\,$ for $\,t\geq0\,$ is a local martingale.
	\end{itemize}
    In particular, Theorem~\ref{thm:Watanabe} implies that $N$ is a conditional Poisson process with intensity $YX$ up to time $S$ under $\mathbb{Q}$.
\begin{proof}
    We show that $\,\mathbb{E}[Z_t]=1\,$ for finite $\,t\geq0$. Then, statements (i) and (ii) follow from well-known results; see, for example, \textcite[Ch VI.2]{Bremaud1981}.

    We show by induction over $\,j\in\mathbb{N}_0\,$ with $\,n\geq j\geq0\,$ and the convention $\,0^0=1\,$ that
    \begin{equation}\label{eq:E_prod=}
        \mathbb{E}\Big[\prod_{k=1}^nY_{T_k}\,\mathbbm{1}_{\{N_t=n\}}\,\Big| \,
        \mathcal{F}_{T_{n-j}}\Big] = \mathrm{e}^{\textstyle-\int_{T_{n-j}}^tX_s\,\mathrm{d}s}
        \frac{1}{j!} \Big(\int_{T_{n-j}}^t Y_sX_s\,\mathrm{d}s\Big)^j 
        \,\mathbbm{1}_{\{t\geq T_{n-j}\}} \, \prod_{k=1}^{n-j}Y_{T_k}.
    \end{equation}
    The base case $\,j=0\,$ follows from \eqref{eq:EW_Tn+1=} in Remark~\ref{rmk:each-jump} for the specific choice $\,W_s=\mathbbm{1}_{\{s>t\geq T_n\}}$ together with the identity $\,\int_{T_n}^\infty\Phi_s^n\,\mathrm{d}s+\Phi_\infty^n=\exp(-\int_{T_n}^t\!X_s\,\mathrm{d}s)$ on $\,\{t\geq T_n\}$. More precisely, we have
    \begin{align*}
    \mathbb{E}\Big[ \prod_{k=1}^n Y_{T_k}\,\mathbbm{1}_{N_t=n}\,\Big|\,\mathcal{F}_{T_n}\Big] 
    &=\mathbb{P}\big[T_{n+1}>t\geq T_n\,\big|\,\mathcal{F}_{T_n}\big] \, \prod_{k=1}^nY_{T_k} \\
    &=\mathrm{e}^{\textstyle-\int_{T_n}^t X_s\,\mathrm{d}s} \, 
    \mathbbm{1}_{\{t\geq T_n\}} \, \prod_{k=1}^nY_{T_k}.
    \end{align*}    
    Let us continue with the induction step; we assume that \eqref{eq:E_prod=} holds for $\,j-1$;
    \begin{align*}
        &\mathbb{E}\Big[\prod_{k=1}^nY_{T_k} \,\mathbbm{1}_{\{N_t=n\}} \,\Big|\,\mathcal{F}_{T_{n-j}}\Big]
        =\mathbb{E}\Big[\mathbb{E}\Big[\prod_{k=1}^nY_{T_k}\,\mathbbm{1}_{\{N_t=n\}}\,\Big|\,\mathcal{F}_{T_{n-(j-1)}} \Big] \,\Big|\,\mathcal{F}_{T_{n-j}}\Big]\\ 
        &\stackrel{\eqref{eq:E_prod=}}{=}
        \mathbb{E}\Big[Y_{T_{n-(j-1)}}\,\mathrm{e}^{\textstyle-\int_{T_{n-(j-1)}}^t\!X_s\,\mathrm{d}s}\frac{1}{(j-1)!} \Big(\int_{T_{n-(j-1)}}^t \!\!\!\! \!\!\!\! \!\! Y_sX_s\,\mathrm{d}s\Big)^{j-1}\;\mathbbm{1}_{\{t\geq T_{n-(j-1)}\}}\,\Big|\,\mathcal{F}_{T_{n-j}}\Big]\,\prod_{k=1}^{n-j}Y_{T_k}\\
        &\stackrel{\eqref{eq:EW_Tn+1=}}{=} \int_{T_{n-j}}^tY_r\,X_r\,\mathrm{e}^{\textstyle-\int_{T_{n-j}}^rX_s\,\mathrm{d}s} \,\mathrm{e}^{\textstyle-\int_r^tX_s\,\mathrm{d}s} \frac{1}{(j-1)!}\Big(\int_r^tY_sX_s\,\mathrm{d}s\Big)^{j-1}\;\mathbbm{1}_{t\geq T_{n-j}}\,\mathrm{d}r\,\prod_{k=1}^{n-j}Y_{T_k}\\
        &=\mathrm{e}^{\textstyle-\int_{]T_{n-j},t]}X_s\,\mathrm{d}s} \, \int_{]{T_{n-j},t}]} \!\!\!\! Y_r\,X_r\,\frac{1}{(j-1)!} \Big( \int_{]r,t]}Y_sX_s\,\mathrm{d}s\Big)^{j-1}\;\mathrm{d}r\;\;\mathbbm{1}_{t\geq T_{n-j}}\,\prod_{k=1}^{n-j}Y_{T_k} \\
        &=\mathrm{e}^{\textstyle-\int_{T_{n-j}}^tX_s\,\mathrm{d}s} \, \frac{1}{j!} \, \Big( \int_{T_{n-j}}^tY_sX_s\,\mathrm{d}s \Big)^j \; \mathbbm{1}_{t\geq T_{n-j}} \, \prod_{k=1}^{n-j}Y_{T_k},
   \end{align*}
    which concludes the induction proof.

    Now, taking $\,j=n\,$ in \eqref{eq:E_prod=} and summing over $n$ yields
    \begin{align*}
        \sum_{n=0}^\infty \mathbb{E}\Big[\prod_{k=1}^nY_{T_k} \,\mathbbm{1}_{N_t=n}\,\Big|\,\mathcal{F}_0\Big]
        &= \mathrm{e}^{\textstyle-\int_0^tX_s\,\mathrm{d}s} \sum_{n=0}^\infty\frac{1}{n!} \,\Big(\int_0^tY_sX_s\,\mathrm{d}s\Big)^n\\
        &=\exp\Big(\int_0^t(Y_s-1)\,X_s\,\mathrm{d}s\Big).
    \end{align*}
    Hence, after applying the tower property and monotone convergence
    \begin{align*}
        \mathbb{E}[Z_t]&=\mathbb{E}
        \Big[\exp\Big(-\int_0^t(Y_s-1)X_s\,\mathrm{d}s\Big) \,\sum_{n=0}^\infty \mathbb{E}\Big[\prod_{k=1}^nY_{T_k} \, \mathbbm{1}_{N_t=n}\,\Big|\,\mathcal{F}_0\Big]\Big]
        \\&=\mathbb{E}\Big[\exp\Big(-\int_0^t(Y_s-1)\,X_s\,\mathrm{d}s\Big)\,
        \exp\Big(\int_0^t(Y_s-1)\,X_s\,\mathrm{d}s\Big)\Big]=1.
    \end{align*}
\end{proof}
\end{thm}

\begin{rmk}\label{rmk:deterministic-variation}
    We can also show $\,\mathbb{E}[Z_t]=1\,$ from the proof of Theorem~\ref{thm:Change-of-Intensity} using general semimartingale theory.% by slightly strengthening the assumptions. 

    Because $Z_t$ is supposed to change the distribution of $N$ in a particular way, we know the distribution of $N$ before and the supposed distribution of $N$ after the change of measure in terms of its semimartingale characteristics, regardless of the soundness of the change of measure. In particular, one can check whether one distribution is absolutely continuous with respect to the other and, if so, whether $Z_t$ coincides with the corresponding Radon–Nikodym derivative. 

    \textcite[Theorem 3.1, Proposition 3.12]{KallsenKarbe2010} used this approach to establish a sufficient condition when a stochastic exponential of a local martingale $M$ is a true martingale in terms of $M$'s semimartingale characteristics. In particular, they showed that the stochastic exponential of $M$ is a true martingale assuming that $M$ has independent increments and fulfils some mild restrictions. It is worth noting that \textcite[Theorem~4.1(3)]{CR} showed that the mild restrictions are unnecessary and the martingale property follows already from the independence of increments. Because we have
    \[
        Z=\mathcal{E}\Big(\int_0^\cdot Y_s-1\,\mathrm{d}\widetilde{N}_s\Big)    
    \] 
    and $\,M_t=\int_0^tY_s-1\,\mathrm{d}\widetilde{N}_s\,$ fulfils the aforementioned assumptions conditioned on $\mathcal{F}_0$, their result yields $\,\mathbb{E}[Z_t\,|\,\mathcal{F}_0]=1$.

    %T: incorporated new citation
    
    %\textcite[Theorem 3.1, Proposition 3.12]{KallsenKarbe2010} used this approach to establish a sufficient condition when a stochastic exponential of a local martingale $M$ is a true martingale in terms of $M$'s semimartingale characteristics. In particular, they showed that the stochastic exponential of $M$ is a true martingale assuming that $M$ has independent increments and jumps with $\,\Delta M_s>-1$. Because we have
    %\[
    %    Z=\mathcal{E}\Big(\int_0^\cdot Y_s-1\,\mathrm{d}\widetilde{N}_s\Big)    
    %\] 
    %and $\,M_t=\int_0^tY_s-1\,\mathrm{d}\widetilde{N}_s\,$ fulfils the aforementioned assumptions conditioned on $\mathcal{F}_0$ under the restriction that $\,Y_s>0$, their result yields $\,\mathbb{E}[Z_t\,|\,\mathcal{F}_0]=1$ assuming $\,Y_s>0$.
\end{rmk}

We finish our note with the following example, which illustrates the reference probability approach and its connection to Theorem~\ref{thm:Change-of-Intensity}.

\begin{ex}[Filtering compound Poisson processes from observations about the jump times]\label{ex:MSc}
    Consider $\,Y=1/X\,$ (assuming $\,X_s>0$). Then, Theorem~\ref{thm:Change-of-Intensity} yields $\mathbb{Q}$ under which $N$ has a constant intensity of $1$ in a given interval $[0,T]$. Hence, $N$ is a Poisson process in $[0,T]$ under $\mathbb{Q}$ according to \textcite[Remark to Theorem 2.3]{Watanabe1964}, so future increments $\,N_t-N_s$, $\,0\le s< t\le T$, are independent of current information $\mathcal{F}_s$ under $\mathbb{Q}$.
    %T: It should be F not F^N

    The Kallianpur--Striebel formula (see \textcite[Proposition 3.16]{BainCrisan2009}) allows us to rewrite any conditional expectation under $\mathbb{P}$ in terms of conditional expectations under $\mathbb{Q}$. More precisely, let $f$ be a positive Borel measurable function and $\,(\mathcal{F}^N_t)_{t\ge0}$ be the natural filtration of $N$ (in particular, not including information about $X$) and note that Theorem~\ref{thm:Change-of-Intensity} implies $\,\mathrm{d}\mathbb{Q}/\mathrm{d}\mathbb{P}|_{\mathcal{F}_t}=Z_t\,$ for $\,T\geq t\geq 0$, then
    \[\mathbb{E}[f(X_t)\,|\,\mathcal{F}_t^N]=\frac{\mathbb{E}_\mathbb{Q}[f(X_t)\,\mathrm{d}\mathbb{P}/\mathrm{d}\mathbb{Q}\,|\,\mathcal{F}_t^N]}{\mathbb{E}_\mathbb{Q}[\mathrm{d}\mathbb{P}/\mathrm{d}\mathbb{Q}\,|\,\mathcal{F}_t^N]}=\frac{\mathbb{E}_\mathbb{Q}[f(X_t)/Z_t\,|\,\mathcal{F}_t^N]}{\mathbb{E}_\mathbb{Q}[1/Z_t\,|\,\mathcal{F}_t^N]}.\] 
    %D: notation $\mathcal{F}^N$ ?: introduced where, why different to p9 here? %T: Good catch; it is introduced on page 9 now and made consistent with the rest of the example 
    Observe that $f(X_t)/Z_t$ and $1/Z_t$ are explicit functions of $N$ and $X$ only. In particular, we can compute the conditional expectations on the right side of the above equation by integrating $X$ while treating $N$ as a known fixed parameter because $X$ and $N$ are independent under $\mathbb{Q}$ and $\,\mathbb{Q}=\mathbb{P}\,$ on $\mathcal{F}_0$. 

    For example, let $X$ be a compound Poisson process starting initially in $\,x_0>0\,$ with an independent Poisson process $M$ and positive jumps $\,(\xi_n)_{n=1}^\infty$, i.e.\ $X_t=x_0+\sum_{k=0}^{M_t}\xi_k\,$ for $\,t\geq0$. Assume that we want to describe the distribution of $X_t$ in terms of its Laplace transform given the observation $\mathcal{F}^N_t$  and additional knowledge $\mathcal{F}^M_t$ (natural filtration of $M$) of the jump times up to time $t$, i.e.\ $\mathbb{E}[\exp(\alpha X_t)\,|\,\mathcal{F}_t^N\vee\mathcal{F}_t^M]\,$ for $\,\alpha\in\mathbb{R}$. Then, it can be shown that
    %D: more precise details. I do not understand, e.g. filtration notations introduced where? raw/natural or usual conds generated? question: Do we require usual conditions in general, see p1. and actually need that filtration is generated by $N$? I doubt so. % T: F^N stands for the natural filtration of N (similar F^M and M) and F is the underlying filtration that includes the knowledge about X (for simplicity, F assumed to be complete - so it does not create issues with stopping times); comments and details added; hope that it is clearer now 
    \begin{gather*}
        \mathbb{E}\big[\exp\big(\alpha X_t\big)\,\big|\,\mathcal{F}_t^N\vee \mathcal{F}_t^M\big]=\frac{\displaystyle\Lambda_t^\alpha}{\Lambda_t^0}\quad\mbox{with}
        \\\Lambda_t^\alpha=\int_{\mathbb{R}^m}\exp\big(\alpha X^j_t\big)\, \exp\bigg(\int_0^t(1-X_s^j)\,\mathrm{d}s\bigg)\prod_{s\in]0,t]}\!\big(1+(X_s^j-1)\,\Delta N_s\big)\;\mathrm{d}\Gamma_{\!m}(j)\bigg|_{m=M_t},
    \end{gather*}
    in which $\,M_t=m\,$ determines the number of integrals in the multiple integral, i.e.\ $j\in\mathbb{R}^m$, $\Gamma_{m}$ is the joint distribution of $(\xi_k)_{k=1}^m$, and $\,X^j_s= x_0+\sum_{k=1}^{M_s}j_k$, i.e.\ $X$ with fixed jump sizes determined by the components of $j$.
    %D: where from / what is $j, j_k, X^j$...? %T: added more details.
\end{ex}

%-----------------------
\printbibliography%biblatex uses different command to insert bibliography
%-----------------------

\end{document}